\documentclass[12pt]{article}

\usepackage{amsfonts}
\usepackage{xr}
\usepackage{graphicx,psfrag}
\usepackage{amsmath}
\usepackage{color,soul,textcomp}
\usepackage[active]{srcltx}

\def\ds{\displaystyle}
\def\2{C^{1,2}(\R\times\R^N)}
\def\to{\rightarrow}
\def\e{\varepsilon}

\def\A{\mathcal{A}}
\def\C{\mathcal{C}}

\def\R{\mathbb{R}}

\def\tilde{\widetilde}
\def\.{\cdot}

\def\lp {\left( }
\def\rp {\right) }

\newlength{\textlarg}

\newcommand{\be}{\begin{equation}}
\newcommand{\ee}{\end{equation}}
\newcommand{\baa}{\begin{array}}
\newcommand{\eaa}{\end{array}}
\newcommand{\ba}{\begin{eqnarray}}
\newcommand{\ea}{\end{eqnarray}}

\newtheorem{thm}{\bf Theorem}[section]
\newtheorem{lem}[thm]{\bf Lemma}
\newtheorem{prop}[thm]{\bf Proposition}

\newtheorem{rmq}[thm]{\bf Remark}

\newenvironment{formula}[1]{\begin{equation}\label{eq:#1}}
                       {\end{equation}\noindent}
\def\Fi#1{\begin{formula}{#1}}
\def\Ff{\end{formula}\noindent}

\begin{document}
\date{}
\title{Maximal and minimal spreading speeds for reaction diffusion equations in nonperiodic slowly varying media}

\author{Jimmy Garnier$^{\hbox{ \small{a,b}}}$, Thomas Giletti$^{\hbox{ \small{b}}}$, Gregoire~Nadin$^{\hbox{ \small{c}}}$\\
\footnotesize{$^{\hbox{a }}$UR 546 Biostatistique et Processus Spatiaux, INRA, F-84000 Avignon, France}\\
\footnotesize{$^{\hbox{b }}$Aix-Marseille Universit\'e, LATP UMR 6632, Facult\'e des Sciences et Techniques}\\
\footnotesize{Avenue Escadrille Normandie-Niemen, F-13397 Marseille Cedex 20, France}\\
\footnotesize{$^{\hbox{c }}$
CNRS, UMR 7598, Laboratoire Jacques-Louis Lions, F-75005 Paris, France}
}

\maketitle

\begin{abstract}
This paper investigates the asymptotic behavior of the solutions of the Fisher-KPP equation in a heterogeneous medium,
$$\partial_t u = \partial_{xx} u + f(x,u),$$
associated with a compactly supported initial datum.
A typical nonlinearity we consider is $f(x,u) = \mu_0 (\phi (x)) u(1-u)$, where $\mu_0$ is a $1$-periodic function and $\phi$ is a $\mathcal{C}^1$ increasing function 
that satisfies $\lim_{x\to +\infty} \phi (x) = +\infty$ and $\lim_{x\to +\infty} \phi' (x) = 0$. 
Although quite specific, the choice of such a reaction term is motivated by its highly heterogeneous nature. We exhibit two different behaviors for $u$ for large times, depending on 
the speed of the convergence of $\phi$ at infinity. If $\phi$ grows sufficiently slowly, then we prove that the spreading speed of $u$ oscillates between two distinct values. 
If $\phi$ grows rapidly, then we compute explicitly a unique and well determined speed of propagation $w_\infty$, arising from the limiting problem of an infinite period.
We give a heuristic interpretation for these two behaviors. 

\end{abstract}


\smallskip


\smallskip

\noindent This work was partially supported by the French ANR project {\em Prefered}.


\section{Introduction}

\subsection{Hypotheses}

We consider the following reaction-diffusion equation in $(0,+\infty) \times \mathbb{R}$:
\begin{equation}\label{eqn:eqRD}
\partial_t u = \partial_{xx} u + f(x,u).
\end{equation}
We assume that $f=f(x,u)$ is locally Lipschitz-continuous in $u$ and of class $\C^1$ in the neighborhood of $u=0$ uniformly with respect to $x$, so that we can define 
$$\mu (x) := f'_u (x,0).$$ Moreover, $f$ is of the KPP type, that is 
$$ f(x,0) =0, \ f(x,1) \leq 0, \ \mu (x) >0 \ \hbox{ and } f(x,u) \leq \mu (x) u \ \hbox{ for all } (x,u) \in \R \times (0,1).$$
A typical $f$ which satisfies these hypotheses is $f(x,u) = \mu (x) u(1-u)$, where $\mu$ is a continuous, positive and bounded function. 

The very specific hypothesis we make on $f$ in this paper is the following: there exist $\mu_0\in \C^0(\R)$ and $\phi\in \C^1 (\R)$ such that
\begin{equation}\label{eqn:hypmu}
\left\{
\begin{array}{l}
\ds \mu (x) = \mu_0 (\phi (x)) \hbox{ for all } x\in\R, \\[0.2cm]
\ds 0 < \min_{[0,1]} \mu_0 < \max_{[0,1]} \mu_0 \  \mbox{ and } \ \mu_0 \mbox{ is 1-periodic}, \vspace{3pt} \\
\ds \phi '(x) > 0, \  \lim_{x\rightarrow +\infty} \phi (x) = +\infty \  \mbox{ and } \lim_{x \rightarrow +\infty} \phi ' (x) =0 .
\end{array}
\right.
\end{equation}

That is, our reaction-diffusion equation is strictly heterogeneous (it is not even almost periodic or ergodic), which means that it can provide useful information on both efficiency of recently developed tools and properties of the general heterogeneous problem. But it also satisfies some periodicity properties with a growing period near $+\infty$.
We aim to look at the influence of the varying period $L(x) := x/\phi (x)$ on the propagation of the solutions.

Note that we do not assume here that there exists a positive stationary solution of~$(\ref{eqn:eqRD})$. We require several assumptions that involve the linearization of $f$
near $u=0$ but our only assumption which is related to the behavior of $ f= f(x,u)$ with respect to $u>0$ is that $f(x,1)\leq 0$, that is, $1$ is a supersolution of~$(\ref{eqn:eqRD})$ (it is clear that, up to some change of variables,
 $1$ could be replaced by any positive constant in this inequality). It is possible to prove that there exists a minimal and stable positive stationary solution of $(\ref{eqn:eqRD})$
by using this hypothesis and the fact that $\mu_0$ is positive~\cite{BHR07}, but we will not discuss this problem since this is not the main topic of this paper.

\subsection{Definitions of the spreading speeds and earlier works}

For any compactly supported initial condition $u_0$ with $0\leq u_0 \leq 1$ and $u_0 \not \equiv 0$, we define the {\em minimal and maximal spreading speeds} as:
\begin{eqnarray*}
w_* & = & \sup \{ c>0 \ | \ \liminf \; \inf_{x \in [0,ct]} u(t,x) > 0 \mbox{ as } t \rightarrow +\infty \}, \vspace{3pt}\\
w^* &= &  \inf \{ c>0 \ | \ \sup_{x \in [ct,+\infty)} u(t,x) \rightarrow 0 \mbox{ as } t \rightarrow +\infty\}.
\end{eqnarray*}
Note that it is clear, from the strong maximum principle, that for any $t >0$ and $x \in \R$, one has $0 < u(t,x) < 1$. One can also easily derive from the homogeneous case~\cite{AronsonWeinberger} that 
$$2\sqrt{\min \mu_0} \leq w_* \leq w^* \leq 2\sqrt{\max \mu_0}.$$
The reader could also remark that we just require $\liminf_{t\to +\infty} u(t,x+ct) > 0 $ in the definition of $w_*$. This is because we did not assume the existence of a 
positive stationary solution. Hence, we just require $u$ to ``take off'' from the unstable steady state $0$. 

The aim of this paper is to determine if some of these inequalities are indeed equalities.\\

\smallskip

The first result on spreading speeds is due to Aronson and Weinberger \cite{AronsonWeinberger}. They proved that 
$w^*=w_* = 2\sqrt{f'(0)}$ in the case where $f$ does not depend on $x$. More generally, even if $f$ does not satisfy $f(u) \leq f'(0) u$ for all $u\in [0,1]$, then 
$w^*=w_*$ is the minimal speed of existence of traveling fronts \cite{AronsonWeinberger}. However, because of the numerous applications in various fields of natural sciences, the role of heterogeneity has become an important topic in the mathematical analysis.

When $f$ is periodic in $x$, Freidlin and Gartner
 \cite{GartnerFreidlin} and Freidlin \cite{Freidlin} 
proved that $w_*= w^*$ using probabilistic techniques. In this case, the spreading speed is 
characterized using periodic principal eigenvalues. Namely, assume that $f$ is $1$-periodic in $x$, set $\mu_0 (x):= f_u'(x,0)$ and define for all~$p\in\R$ the elliptic operator
\begin{equation} \mathcal{L}_p \varphi := \varphi'' - 2p\varphi ' + (p^2 + \mu_0 (x))\varphi.\end{equation}
It is known from the Krein-Rutman theory that this operator admits a unique periodic principal eigenvalue $\lambda_p (\mu_0)$, defined by the existence of a positive $1$-periodic 
function $\varphi_p\in\mathcal{C}^2 (\R)$ so that $\mathcal{L}_p \varphi_p = \lambda_p (\mu_0)\varphi_p$. The characterization of the spreading speed \cite{GartnerFreidlin} reads 
\begin{equation} \label{eq:characterizationw} w_* = w^* =\min_{p>0} \frac{\lambda_p (\mu_0)}{p}.\end{equation}
Such a formula is very useful to investigate the dependence between the spreading speed and the growth rate $\mu_0$. Several alternative proofs of this characterization, 
based on different techniques, have been given in \cite{BHNa, Weinberger}. The spreading speed $w_*=w^*$ has also been identified later
as the minimal speed of existence of pulsating traveling fronts, which is the appropriate generalization of the notion of traveling fronts to periodic media 
\cite{BerestyckiHamel}. Let us mention, without getting into details, that the equality $w_*=w^*$ and the characterization~(\ref{eq:characterizationw}) have been extended 
when the heterogeneity is transverse~\cite{MallordyRoquejoffre}, space-time periodic or compactly supported~\cite{BHNa}, or random stationary ergodic 
\cite{GartnerFreidlin, NolenXinrandom}. In this last case one has to use Lyapounov exponents instead of principal eigenvalues. 

In all these cases (except in the random one), the operator $\mathcal{L}_p$ is compact and thus principal eigenvalues are well-defined. When the dependence of $f$ with respect to $x$
is more general, then classical principal eigenvalues are not always defined, which makes the computation of the spreading speeds much more difficult.
Moreover, in general heterogeneous media, it may happen that $w_*<w^*$. No example of such phenomenon has been given in space heterogeneous media, but there exist 
examples in time heterogeneous media \cite{BerestyckiNadin} or when the initial datum is not compactly supported \cite{HamelNadin}. 

Spreading properties in general heterogeneous media have recently been investigated by Berestycki, Hamel and the third author in \cite{BHNa}. These authors 
clarified the links between the different notions of spreading speeds and gave some estimates on the spreading speeds. 
More recently, Berestycki and the third author gave sharper bounds using the notion of generalized principal eigenvalues \cite{BerestyckiNadin}. 
These estimates are optimal when the nonlinearity is periodic, almost periodic or random stationary ergodic. In these cases, one gets
$w_*=w^*$ and this spreading speed can be characterized through a formula which is similar to (\ref{eq:characterizationw}), 
involving generalized principal eigenvalues instead of periodic principal eigenvalues.


\section{Statement of the results}

Before enouncing our results, let us first roughly describe the situation. As $\phi'(x)\to 0$ as $x\to +\infty$, the function $\phi$ is sublinear at infinity and thus $\mu (x) = \mu_0 (\phi (x))$ stays near its extremal values 
$\max \mu_0$ or $\min \mu_0$ on larger and larger intervals. If these intervals are sufficiently large, that is, if $\phi$ increases sufficiently slowly, 
the solution $u$ of (\ref{eqn:eqRD}) should propagate alternately at speeds close to $2\sqrt{\max \mu_0}$ and $2\sqrt{\min \mu_0}$.
Hence, we expect in such a case that $w^*=2\sqrt{\max \mu_0}$ and $w_*=2\sqrt{\min \mu_0}$. 
%

On the other hand, if one writes $\phi (x) = x/ L(x)$, then the reaction-term locally looks like an $L(x)$-periodic function. Since $L(x) \to +\infty$, as clearly follows from the fact that $\phi'(x) \to0$ as $x\to +\infty$,
one might expect to find a link between the spreading speeds and the limit of the spreading speed $w_L$ associated with the $L$-periodic 
growth rate $\mu_L (x):= \mu_0 (x/L)$ when $L\to +\infty$. 
This limit has recently been computed by Hamel, Roques and the third author~\cite{HamelNadinRoques}. 
As $\mu_L$ is periodic, $w_L$ is characterized by (\ref{eq:characterizationw}) and one can compute the limit of~$w_L$ by computing the limit of $\lambda_p (\mu_L)$ for 
all $p$. This is how the authors of \cite{HamelNadinRoques} proved that 
\begin{equation}
\lim_{L\to +\infty} w_L = \min_{k \geq M}\frac{k}{j(k)},
\end{equation}
where $M := \max_{x\in\R} \mu_0 (x) >0 $ and $j: [M,+\infty) \to [j(M),+\infty)$ is defined for all $k\geq M$ by 
\begin{equation} \label{eq:defj} j(k):= \int_0^1 \sqrt{k- \mu_0 (x)}dx.\end{equation}
If $\phi$ increases rapidly, that is, the period $L(x)$ increases slowly, then we expect to recover this type of behavior. More precisely, we expect that $w^*= w_* = \min_{k \geq M}k/j(k)$. 

We are now in position to state our results. 

\subsection{Slowly increasing $\phi$}

We first consider the case when $\phi$ converges very slowly to $+\infty$ as $x \rightarrow +\infty$. As expected, we prove in this case that 
$w_* < w^*$.

\begin{thm}\label{th:gencase} 
\begin{enumerate}
\item Assume that $\ds\frac{1}{x \phi '(x)} \rightarrow +\infty$ as $x \rightarrow +\infty$. Then $$w_*= 2\sqrt{\min \mu_0} < w^* = 2\sqrt{\max \mu_0} .$$
\item Assume that $\ds\frac{1}{x \phi '(x)}\to C$ as $x \rightarrow +\infty$. If $C$ is large enough (depending on $\mu_0$), then $$ w_* < w^* .$$
\end{enumerate}
\end{thm} 

This is the first example, as far as we know, of a space heterogeneous nonlinearity $f(x,u)$ for which the spreading speeds $w_*$ and $w^*$ 
associated with compactly supported initial data are not equal.

In order to prove this Theorem, we will first consider 
the particular case when $\mu_0$ is discontinuous and only takes two values (see Proposition \ref{th:2values} below). In this case, we are able to construct sub- and super-solutions on each
interval where $\mu$ is constant, and to conclude under some hypotheses on the length of those intervals. Then, in the general continuous case, our hypotheses on 
\textbf{$\ds(x \phi '(x))^{-1}$} allow us to bound $\mu$ from below and above by some two values functions, and our results then follow from the preliminary case.

\begin{rmq}Note that such a two values case is not continuous, so that our Theorem holds under more general hypotheses. In fact, one would only need that $\mu_0$ is continuous on two points such that $\mu_0$ attains its maximum and minimum there, so that, from the asymptotics of $\phi (x)$, the function $\mu (x)=\mu_0 (\phi(x))$ will be close to its maximum and minimum on very large intervals as $x \rightarrow +\infty$.
\end{rmq}

\subsection{Rapidly increasing $\phi$}

We remind the reader that $M := \max_{x\in\R} \mu_0 (x) >0 $ and $j: [M,+\infty) \to [j(M),+\infty)$ is defined by (\ref{eq:defj}). We expect to characterize the spreading speeds 
$w_*$ and $w^*$ using these quantities, as in \cite{HamelNadinRoques}. 

Note that $j (M)> 0$ since $\min \mu_0 < M$. The function $j$ is clearly a bijection and thus one can define
\begin{equation}
 w_\infty := \min_{\lambda \geq j(M)}\frac{j^{-1}(\lambda)}{\lambda}=\min_{k \geq M}\frac{k}{j(k)}. 
\end{equation}
We need in this section an additional mild hypothesis on $f$:
\begin{equation} \label{hyp:fHolder} \exists\, C>0, \gamma>0 \hbox{ such that } \ f(x,u)\geq f_u'(x,0)u -Cu^{1+\gamma} \  \hbox{ for all } (x,u)\in \R\times(0,+\infty).\end{equation}

\begin{thm} \label{thm-vitspread}
Under the additional assumptions $(\ref{hyp:fHolder})$, $\phi\in \C^3 (\R)$ and 
\begin{equation}\label{hyp:phi2}
\phi'' (x)/\phi' (x)^2 \to 0, \ \hbox{ and}\  \phi''' (x)/\phi' (x)^2 \to 0, \hbox{ as } x\to+\infty , 
\end{equation}
one has
$$w_*= w^* = w_\infty.$$
\end{thm}
Note that (\ref{hyp:phi2}) implies $\ds(x \phi '(x))^{-1}\to 0$ as $x\to +\infty$. Hence, this result is somehow complementary to Theorem \ref{th:gencase}. However, this is not optimal as this does not cover all cases. An interesting and open question would be to refine those results to get more precise necessary and sufficient conditions for the equality $w_* = w^*$ to be satisfied. This could provide some insight on the general heterogeneous case, where the establishment of such criteria is an important issue.

This result will mainly be derived from Theorem 2.1 of \cite{BerestyckiNadin}. We first construct some appropriate test-functions using the asymptotic problem associated with
$\mu_L (x) = \mu_0 (x/L)$ as $L\to +\infty.$ This will enable us to compute the generalized principal eigenvalues and the 
computation of the spreading speeds will follow from 
Theorem 2.1 in \cite{BerestyckiNadin}.

\subsection{Examples}

We end the statement of our results with some examples which illustrate the different possible behaviors. 

\bigskip

\noindent {\bf Example 1:} $\phi (x)=\beta (\ln x)^\alpha,$ with $\alpha, \beta >0$. This function clearly satisfies the hypotheses in (\ref{eqn:hypmu}). 
\begin{itemize}
\item If $\alpha \in (0,1)$, one has 
$1/(x\phi'(x)) = (\ln x)^{1-\alpha} / (\beta \alpha) \to +\infty$ as $x\to +\infty$.
Hence, the assumptions of case $1$ in Theorem \ref{th:gencase} are satisfied and one has $w_*= 2\sqrt{\min\mu_0}$ and $w^* = 2\sqrt{\max \mu_0}.$
\item If $\alpha =1$, then $x\phi'(x) = \beta$ for all $x$ and thus we are in the framework of case $2$ in Theorem \ref{th:gencase}, 
which means that we can conclude that $w_* < w^*$ provided that $\beta$ is small enough. 
\item Lastly, if $\alpha >1$, then straightforward computations give
$$\phi''(x)/\phi'(x)^2 \sim -\frac{1}{\beta \alpha} (\ln x)^{1-\alpha} \to 0 \hbox{ as } x\to +\infty,$$
$$\phi'''(x)/\phi'(x)^2 \sim \frac{2}{\beta \alpha x} (\ln x)^{1-\alpha} \to 0 \hbox{ as } x\to +\infty.$$
Hence, the assumptions of Theorem \ref{thm-vitspread} are satisfied and there exists a unique spreading speed: $w_*= w^* = w_\infty.$

\end{itemize}

\bigskip

\noindent {\bf Example 2: $\phi (x)=x^\alpha, \alpha \in (0,1)$.} This function clearly satisfies the hypotheses in (\ref{eqn:hypmu}) since $\alpha<1$. One has 
$\phi'' (x)/\phi'(x)^2 = \frac{\alpha -1}{\alpha x^\alpha} \to 0$ and $\phi''' (x)/\phi'(x)^2 = \frac{(\alpha -1)(\alpha -2)}{\alpha x^{1+\alpha}} \to 0$ as $x\to +\infty$.
Thus, the assumptions of Theorem \ref{thm-vitspread} are satisfied and $w_*= w^* = w_\infty.$

\bigskip

\noindent {\bf Example 3: $\phi (x)=x/(\ln x)^\alpha, \alpha >0$.} This function satisfies (\ref{eqn:hypmu}) and one has
$$\phi '(x) = \frac{1}{(\ln x)^\alpha} -\frac{\alpha } {(\ln x)^{\alpha+1}},$$
$$\phi ''(x) = \frac{-\alpha}{x(\ln x)^{1+\alpha}} +\frac{\alpha (\alpha+1) } {x(\ln x)^{\alpha+2}},$$
$$\phi '''(x) = \frac{\alpha}{x^2(\ln x)^{1+\alpha}} -\frac{\alpha (\alpha+1)(\alpha +2) } {x^2 (\ln x)^{\alpha+3}}.$$
It follows that $\phi''(x)/\phi'(x)^2 \to 0$ and  $\phi'''(x)/\phi'(x)^2 \to 0$ as $x\to +\infty$ since the terms in $x$ will decrease faster than the terms in $\ln x$. 
Thus, the assumptions of Theorem \ref{thm-vitspread} are satisfied and $w_*= w^* = w_\infty.$

\bigskip

\noindent {\bf Organization of the paper:} Theorem \ref{th:gencase} will be proved in Section \ref{sec:cont}. As a first step to prove this Theorem, we will 
investigate in Section \ref{sec:twovalues} the case where $\mu_0$ is not continuous anymore but only takes two values $\mu_+$ and $\mu_-$. Lastly, Section \ref{sec:unique}
is dedicated to the proof of Theorem \ref{thm-vitspread}. 

\bigskip

\noindent {\bf Acknowledgements:} The authors would like to thank Fran\c cois Hamel and Lionel Roques for having drawn their attention to the problems investigated in 
this paper. This article was completed while the third author 
was visiting the Department of mathematical sciences of Bath whose hospitality is gratefully acknowledged.


\section{The two values case}\label{sec:twovalues}

We assume first that $\mu$ is discontinuous and only takes two distinct values $\mu_- ,$~$\mu_+ \in (0,+\infty)$. 
Moreover, we assume that 
there exist two increasing sequences $(x_n)_n $ and $(y_n)_n$ such that $x_{n+1} \geq y_n \geq x_n$ for all~$n$, $\lim_{n\to +\infty} x_n = +\infty$ and 
\begin{equation}\label{eqn:hypseq}
\mu (x) = \left\{ \begin{array}{l}
                \mu_+  \mbox{ if } x \in (x_n,y_n) ,\vspace{3pt}\\
                \mu_-  \mbox{ if } x \in (y_n,x_{n+1}).
                \end{array}\right.
\end{equation}
\begin{prop}\label{th:2values} We have:
\begin{enumerate}
\item If $y_n/x_n \rightarrow +\infty$, then $w^*= 2\sqrt{\mu_+}$.
\item If $x_{n+1}/y_n \rightarrow +\infty$, then $w_*= 2\sqrt{\mu_-}$.
\item If $y_n/x_n \rightarrow K>1$, then $w^* \geq \ds 2\sqrt{\mu_+}\frac{K}{(K-1)+ \sqrt{\mu_+ / \mu_-}}$.
\item If $x_{n+1}/y_n \rightarrow K>1$, then $w_* \leq 2\sqrt{\mu_-}\ds\frac{K+\sqrt{\mu_+/\mu_-}}{K+\sqrt{\mu_-/\mu_+}}$.
\end{enumerate}
\end{prop}
It is clear in part 3 (resp. 4) of Proposition~\ref{th:2values} that the lower bound on $w^*$ (resp. upper bound on $w_*$) goes to $2\sqrt{\mu_+}$ (resp. $2\sqrt{\mu_-}$) as $K \rightarrow +\infty.$ Hence, for $K$ large enough, we get the wanted result $w_* < w^*$. 

\subsection{Maximal speed: proof of parts 1 and 3 of Proposition \ref{th:2values}}

1. We first look for a subsolution of equation (\ref{eqn:eqRD}) going at some speed $c$ close to $2\sqrt{\mu_-}$. Let~$\phi_R$ be a solution of the principal eigenvalue problem:
\begin{equation}\label{eq-eigenphiR}
\left\{
\begin{array}{ll}
\displaystyle \partial_{xx} \phi_R  = \lambda_{R} \phi_R & \mbox{ in } B_R,\vspace{3pt}\\
\displaystyle \phi_R=0 & \mbox{ on } \partial B_R,\vspace{3pt}\\
\displaystyle \phi_R>0 & \mbox{ in } B_R.
\end{array}
\right.
\end{equation}
We normalize $\phi_R$ by $\|\phi_R\|_\infty=1$.
We know that $\lambda_{R} \rightarrow 0$ as $R \rightarrow +\infty$. Let $c < 2\sqrt{\mu_-}$ and~$R$ large enough so that
$-\lambda_R < \mu_- - c^2/4$.
Then $\ds v(x)=e^{\frac{-cx}{2}} \phi_R (x)$ satisfies:
$$\partial_{xx} v + c \partial_x v + \mu_- v = \Big( \mu_- - \frac{c^2}{4} + \lambda_R \Big) v >0 \hbox{ in } B_R.$$
By extending $\phi_R$ by~0 outside $B_R$, by regularity of $f$ and since $f'_u (x,0) \geq \mu_-$ for any $x\in \mathbb{R}$, for some small $\kappa$, we also have in $(0,+\infty)\times \mathbb{R}$:
$$\partial_{xx} \kappa v + c \partial_x \kappa v + f(x+ct, \kappa v) \geq 0.$$
Hence, $w (t,x) :=\kappa v (x-ct)$ is a subsolution of (\ref{eqn:eqRD}).
Without loss of generality, we can assume that $u(1,x) \geq w (1,x)$, thus for any $t\geq 1$, $u(t,x) \geq w (t,x)$.
That is, for any speed $c < 2\sqrt{\mu_-}$, we have bounded $u$ from below by a subsolution of (\ref{eqn:eqRD}) with speed c. 
In particular, 
$$\hbox{ let } \ds t_n := \frac{x_n +R}{c}, \hbox{ then } u(t_n,x) \geq w(t_n,x) \hbox{ for all } x\in\R,$$ 
which is positive on a ball of radius $R$ around $x_n +R$.\\

\smallskip

\noindent 2. Take an arbitrary $c'<2\sqrt{\mu_+}$ and let $\phi_{R'}$ a solution of the principal eigenvalue problem~(\ref{eq-eigenphiR}) with $R'$ such that
$-\lambda_{R'} < \mu_+ - c'^2/4$. 
As above, there exists $\ds \tilde{v} (x)=\kappa' e^{\frac{-c' x}{2}} \phi_{R'} (x)$ compactly supported such that
\begin{equation}
\partial_{xx} \tilde{v} + c' \partial_x \tilde{v} + f(x+x_n+R+c't, \tilde{v}) \geq 0, 
\end{equation}
as long as $\tilde{v}=0$ where $f'_u (x+x_n+R+c't,0) \neq \mu_+$, that is $$(-R'+x_n+R+c't,R'+x_n+R+c't) \subset (x_n,y_n),$$ which is true for
$R>R'$ and
$$0\leq t \leq \frac{ y_n -x_n - R -R'}{c'}.$$
As $R$ could be chosen arbitrarily large, we can assume that the condition $R>R'$ is indeed satisfied. 
Moreover, as $\liminf_{n\to +\infty} y_n/x_n >1$ and $\lim_{n \to +\infty} x_n = +\infty$, we can assume that~$n$ is large enough so that $y_n-x_n >2R$ and thus the 
second condition is also satisfied. 
Hence, $\tilde{w} (t_n+t,x):=\tilde{v} (x-x_n-R -c't)$
is a subsolution of (\ref{eqn:eqRD}) for $t~\in~(0,\frac{ y_n -x_n - R -R'}{c'})$ and $x\in \R$.
We can take $\kappa'$ small enough so that
\begin{equation}
\kappa \min_{y\in B(0,R')}\phi_R (y)>\kappa' e^{\frac{|c-c'|}{2}R'}. 
\end{equation}
For all $x\in B(ct_n,R')$, one has:
\begin{equation}
\begin{array}{rcl} w(t_n,x)=\kappa e^{\frac{-c(x-ct_n)}{2}} \phi_R (x-ct_n)
&\geq &\ds \kappa \Big( \min_{y\in B(0,R')}\phi_R (y) \Big)\ e^{\frac{-(c-c')}{2}(x-ct_n)}e^{\frac{-c'}{2}(x-ct_n)}\\[0.2cm]
&\geq &\ds \kappa \Big( \min_{y\in B(0,R')}\phi_R (y)\Big) \ e^{\frac{-|c-c'|}{2}R'}e^{\frac{-c'}{2}(x-ct_n)} \\[0.2cm]
&\geq &\ds \kappa'  e^{\frac{-c'}{2}(x-ct_n)} \\[0.2cm]
&\geq &\ds \kappa' e^{\frac{-c' (x-ct_n)}{2}} \phi_{R'} (x-ct_n)=\tilde{w} (t_n,x), \end{array} 
\end{equation}
since $ct_n = x_n+R$ by definition. 
Moreover, $u(t_n,x)\geq w(t_n,x)$ for all $x\in\mathbb{R}$. The parabolic maximum principle thus gives
$$u(t_n+t,x)\geq \tilde{w} (t_n+t,x) \ \hbox{ for all } t \in \Big(0,\frac{ y_n -x_n - R -R'}{c'}\Big) \hbox{ and } \ x\in\R.$$

\smallskip

  \noindent 3. We can now conclude. Indeed, for $n$ large enough one has:
$$u\left(t_n+\frac{ y_n -x_n - R -R'}{c'},y_n -R'\right) \geq \tilde{w} \left(t_n+\frac{ y_n -x_n - R -R'}{c'},y_n-R'\right)= \tilde{v} (0).$$
Since the construction of $\tilde{v}$ did not depend on $n,$ the above inequality holds independently of $n,$ which implies that:
$$\inf_n u\left(t_n+\frac{ y_n -x_n - R -R'}{c'},y_n -R'\right) > 0.$$ 

If $y_n/x_n \rightarrow +\infty$, we have
$$\frac{y_n - R'}{t_n+\frac{ y_n -x_n - R -R'}{c'}}=\frac{y_n - R'}{\frac{x_n}{c}+\frac{ y_n -x_n - R -R'}{c'}} \rightarrow c' \ \hbox{ as } n \rightarrow +\infty.$$
 It follows that $w^* \geq c'$ for any $c' < 2\sqrt{\mu_+}$. The proof of part 1 of Proposition~\ref{th:2values} is completed.

If $y_n/x_n \rightarrow K$, we have
$$\frac{y_n - R'}{t_n+\frac{ y_n -x_n - R -R'}{c'}} \rightarrow \frac{ K}{\frac{1}{c}+\frac{K-1}{c'}} \ \hbox{ as } n \rightarrow +\infty.$$
As this is true for any $c' < 2\sqrt{\mu_+}$ and $c < 2\sqrt{\mu_-}$, this concludes the proof of part 3 of Proposition \ref{th:2values}.
\hfill$\Box$

\subsection{Minimal speed: proof of parts 2 and 4 of Proposition \ref{th:2values}}

Let $\lambda_+ = \sqrt{\mu_+}$ be the solution of $\lambda_+^2 - 2\sqrt{\mu_+} \lambda_+ = - \mu_+$.
One can then easily check, from the KPP hypothesis, that the function
$$v(t,x):=\min \left(1, \kappa e^{-\lambda_+ (x-2\sqrt{\mu_+}t)} \right)$$
is a supersolution of equation (\ref{eqn:eqRD}) going at the speed $2\sqrt{\mu_+}$, for any $\kappa >0$. Since $u_0$ is compactly supported,
we can choose $\kappa$ such that $v(0,\cdot)\geq u_0$ in $\mathbb{R}$. Thus, for any $t\geq 0$ and $x \in \mathbb{R}$, $u(t,x) \leq v(t,x)$.
In particular, the inequality holds for $t=t_n$ the smallest time such that $v(t,y_n)=1$. Note that $t_n=y_n/(2\sqrt{\mu_+})+C$ where $C$ is a 
constant independent of $n$. 
Then for all $x\in \mathbb{R}$,
$$u(t_n,y_n+x) \leq v(t_n,y_n+x)= \min \left(1, e^{-\lambda_+ x} \right).$$
We now look for a supersolution moving with speed $2\sqrt{\mu_-}$ locally in time around $t_n$. Let us define
$$w(t_n + t,y_n +x):=\min \left(v(t_n +t,y_n+x), e^{-\lambda_- (x-2\sqrt{\mu_-}t)} \right)$$
where $\lambda_- = \sqrt{\mu_-}$. Note that $\lambda_- < \lambda_+$, thus $u(t_n,y_n+x) \leq v(t_n,y_n+x)= w(t_n,y_n+x)$.

We now check that $w$ is indeed a supersolution of equation (\ref{eqn:eqRD}). 
We already know that $v$ is a supersolution and it can easily be seen as above from the KPP hypothesis that $(t,x)\mapsto e^{-\lambda_- (x-2\sqrt{\mu_-}t)}$ is a supersolution only where $f'_u (\cdot,0) =\mu_-$. Thus, we want the inequality
$v(t_n + t,y_n+x) \leq  e^{-\lambda_- (x-2\sqrt{\mu_-}t)}$ to be satisfied if $ y_n +x \not\in (y_n,x_{n+1}).$ Recall that 
$v(t_n+t,y_n+x) = \min \left(1, e^{-\lambda_+ (x-2\sqrt{\mu_+}t)} \right)$ for all $t>0$ and $x\in\R.$
Thus, the inequality is satisfied if $t\geq 0$ and $x\leq 0$ or if
$$x \geq 2 \frac{\lambda_+ \sqrt{\mu_+} - \lambda_- \sqrt{\mu_-}}{\lambda_+ - \lambda_-} t =2(\lambda_++\lambda_-)t.$$
It follows that $w(t_n+t,y_n+x)$ is indeed a supersolution of equation (\ref{eqn:eqRD}) in $\mathbb{R}$ as long as
\begin{equation}\label{eq:tsursol}
0 \leq 2(\lambda_++\lambda_-)t \leq x_{n+1} - y_n,\end{equation}
and that $u(t_n+t,y_n+x) \leq w(t_n+t,y_n+x)$ for any $t$ verifying the above inequality.

To conclude, let now $2\sqrt{\mu_+}>c > 2\sqrt{\mu_-}$, and $t'_n$ the largest $t$ satisfying (\ref{eq:tsursol}), i.e.
$$t_n'= \frac{x_{n+1}-y_n}{2(\lambda_++\lambda_-)}.$$
The sequence $(t_n')_n$ tends to~$+\infty$ as~$n \rightarrow +\infty$ since $\liminf_{n\to +\infty} x_{n+1}/y_n >1$ and $\lim_{n \to +\infty} y_n = +\infty$. 
Moreover, one has 
$$u(t_n+t'_n,y_n+ct'_n) \leq w(t_n+t'_n,y_n+ct'_n) \rightarrow 0\hbox{ as } n \rightarrow +\infty$$
since $c>2\sqrt{\mu_-}$. 

If $x_{n+1}/y_n \rightarrow +\infty$ as $n\to +\infty$, as $t_n=y_n/(2\sqrt{\mu_+})+C$, one gets
$t_n'/t_n\to +\infty$ as~$n\to +\infty$.
Hence,
$$\frac{y_n + ct'_n}{t_n+t'_n} \rightarrow c \hbox{ as } n\to +\infty.$$
It follows that $w_* \leq c$ for any $c > 2\sqrt{\mu_-}$. This proves part 2 of Proposition \ref{th:2values}.

If $x_{n+1}/y_n \rightarrow K$ as $n\to +\infty$, we compute
$$\frac{y_n + ct'_n}{t_n+t'_n} \rightarrow \ds\frac{1+ \frac{c}{2(\lambda_++\lambda_-)}(K-1)}{\frac{1}{2\sqrt{\mu_+}}+\frac{K-1}{2(\lambda_++\lambda_-)}}\hbox{ as } n\to +\infty.$$
Hence, $w_*$ is smaller than the right hand-side. 
As $c\in (2\sqrt{\mu_-},2\sqrt{\mu_+})$ is arbitrary, $\lambda_-= \sqrt{\mu_-}$ and $\lambda_+= \sqrt{\mu_+}$, we eventually get
$$w_* \leq 2\sqrt{\mu_-}\ds\frac{K+\sqrt{\mu_+/\mu_-}}{K+\sqrt{\mu_-/\mu_+}},$$
which concludes the proof of part 4 of Proposition \ref{th:2values}.
\hfill$\Box$


\section{The continuous case}\label{sec:cont}

\subsection{Proof of part 1 of Theorem \ref{th:gencase}}

We assume that $\mu_0$ is a continuous and 1-periodic function. Let now $\varepsilon$ be a small positive constant and define $\mu_- < \mu_+$ by: 
$$\left\{
\begin{array}{rcl}
\mu_+ & := & \max \mu_0 - \varepsilon,\vspace{3pt} \\
\mu_- & := & \min \mu_0.\\
\end{array}
\right.
$$
We want to bound $\mu$ from below by a function taking only the values $\mu_-$ and $\mu_+$, in order to apply Theorem \ref{th:2values}. 
Note first that there exist $x_{-1} \in (0,1)$ and $\delta \in (0,1)$ such that $\mu_0 (x) > \mu_+ $ for any $x \in (x_{-1},x_{-1}+\delta)$. 
We now let the two sequences $(x_n)_{n \in \mathbb{N}}$ and $(y_n)_{n\in \mathbb{N}}$ defined for any $n$ by:
$$\left\{
\begin{array}{rcl}
\phi (x_n) & = & x_{-1} + n,\vspace{3pt}\\
\phi (y_n) & = & x_{-1} + n + \delta.\\
\end{array}
\right.
$$
Note that since $\phi$ is strictly increasing and $\phi(+\infty)=+\infty$, 
then those sequences indeed exist, tend to $+\infty$ as $n \rightarrow +\infty$, and satisfy for any $n$, $x_n < y_n < x_{n+1}$. 
It also immediately follows from their definition that for all $x\in\R,$
$$
\mu (x) \geq \tilde\mu(x) \  \hbox{ where  }\tilde\mu(x):= \left\{ \begin{array}{l}
                \mu_+  \mbox{ if } x \in (x_n,y_n) ,\vspace{3pt}\\
                \mu_-  \mbox{ if } x \in (y_n,x_{n+1}) .\\
                \end{array}\right.
$$
We now have to estimate the ratio $y_n/x_n$ in order to apply Proposition~\ref{th:2values}. Note that:
\begin{equation}\label{eq:div1}
\delta = \phi (y_n) - \phi (x_n) = \int_{x_n}^{y_n} \phi ' (x) dx.
\end{equation}
Moreover, under the hypothesis $x \phi '(x) \rightarrow 0$ as $x \rightarrow +\infty$, and since $(x_n)_n$, $(y_n)_n$ tend to~$+\infty$ as $n \rightarrow +\infty$:
\begin{eqnarray}\label{eq:div2}
\int_{x_n}^{y_n} \phi ' (x) dx  =  \int_{x_n}^{y_n} \left( x \phi ' (x) \ \times \ \frac{1}{x} \right)dx= o \left( \ln \left(\frac{y_n}{x_n}\right) \right) \mbox{ as } n \rightarrow +\infty .
\end{eqnarray}
From (\ref{eq:div1}) and (\ref{eq:div2}), we have that $\frac{y_n}{x_n} \rightarrow + \infty$. To conclude, we use the parabolic maximum principle and part 1 
of Proposition~\ref{th:2values} applied to problem (\ref{eqn:eqRD}) with a reaction term $\tilde{f} \leq f$ such that
$$\tilde{f}'_u (x,0) = \tilde\mu(x) \ \hbox{ for all }x\in\R.
$$
It immediately follows that $w^* \geq 2 \sqrt{\max \mu_0 - \varepsilon}$. Since this inequality holds for any $\varepsilon >0$, we get $w^* = 2 \sqrt{\max \mu_0}$.

We omit the details of the proof of $w_* = 2 \sqrt{\min \mu_0}$ since it follows from the same method.
Indeed, one only have to choose $y_{-1}'$ and $\delta '$ in $(0,1)$ such that $\mu_0 (x) < \min \mu_0 + \varepsilon$ for any $x \in (y_{-1}',y_{-1}'+\delta ')$ 
and let two sequences such that
$$\left\{
\begin{array}{rcl}
\phi (y_n') & = & y_{-1}' + n, \vspace{3pt}\\
\phi (x_{n+1}') & = &y_{-1}' + n + \delta '.\\
\end{array}
\right.
$$
One can then easily conclude as above using part 2 of Proposition~\ref{th:2values}.
\hfill$\Box$


\subsection{Proof of part 2 of Theorem \ref{th:gencase}}

As before, we bound $\mu_0$ from below by a two values function, that is, for all $x \in \R$,
$$
\mu (x) \geq \tilde\mu(x) \  \hbox{ where  }\tilde\mu(x):= \left\{ \begin{array}{l}
                \mu_+ =\max \mu_0 - \e \; \mbox{ if } x \in (x_n,y_n) ,\vspace{3pt}\\
                \mu_-  = \min \mu_0 \; \mbox{ if } x \in (y_n,x_{n+1}) ,\\
                \end{array}\right.
$$
where $\e$ a small positive constant and the two sequences $(x_n)_n$ and $(y_n)_n$ satisfy for any $n$:
$$\left\{
\begin{array}{l}
x_n < y_n < x_{n+1},\vspace{3pt}\\
\phi (x_n) =  x_{-1} + n,\vspace{3pt}\\
\phi (y_n)  =  x_{-1} + n + \delta (\e) \mbox{ for some } \delta (\e) >0,\vspace{3pt}\\
x_n \rightarrow +\infty \mbox{ and } y_n \rightarrow +\infty.
\end{array}
\right.
$$
Here, under the assumption that $x\phi ' (x) \rightarrow 1/C$, we get
\begin{eqnarray*}
\delta (\e) = \phi (y_n) -\phi (x_n) & =&  \int_{x_n}^{y_n} \phi' (x)dx\\
& = & \int_{x_n}^{y_n} \left( x \phi ' (x) \ \times \ \frac{1}{x} \right)dx \\
& = &  \frac{1}{C} \ln \left(\frac{y_n}{x_n}\right) + o \left( \ln \left(\frac{y_n}{x_n}\right) \right)\mbox{ as } n \rightarrow +\infty .
\end{eqnarray*}
Hence,
$$\frac{y_n}{x_n} \rightarrow e^{\delta (\e) C} \mbox{ as } n \rightarrow +\infty.$$
We can now apply the parabolic maximum principle and part~3 of Proposition~\ref{th:2values} to get
\begin{equation}\label{C0-above}
w^* \geq 2\sqrt{\max \mu_0 - \e}\; \frac{e^{\delta (\e) C}}{(e^{\delta (\e) C}-1)+ \sqrt{(\max \mu_0 - \e) / \min \mu_0}}.
\end{equation}
Notice that the dependence of $\delta$ on $\e$ prevents us from passing to the limit as $\e \rightarrow 0$ as we did to prove part~1 of Theorem~\ref{th:gencase}. However, for any fixed $\e >0$, one can easily check that the right-hand side in the inequation (\ref{C0-above}) converges as $C\rightarrow +\infty$ to $2 \sqrt{\max \mu_0 - \e}$.

One can proceed similarly to get an upper bound on $w_*$, that is:
\begin{equation}\label{C0-below}
w_* \leq 2 \sqrt{\min \mu_0 + \e}\; \frac{e^{\delta' (\e) C} + \sqrt{\max \mu_0 / (\min \mu_0 + \e)}}{e^{\delta' (\e) C} + \sqrt{(\min \mu_0 + \e)/\max \mu_0}},
\end{equation}
where $\e$ can be chosen arbitrary small and $\delta ' (\e)$ is such that $\mu_0 (x) \leq \min \mu_0  + \e$ on some interval of length $\delta' (\e)$. It is clear that the right-hand side of (\ref{C0-below}) converges to $2 \sqrt{\min \mu_0 + \e}$ as $C \rightarrow +\infty$.

Therefore, by choosing $\e < (\max \mu_0 - \min \mu_0) /2$, one easily gets from (\ref{C0-above}) and (\ref{C0-below}) that for $C$ large enough, $w_* < w^*$. This concludes the proof of part~2 of Theorem~\ref{th:gencase}. Moreover, note that the choice of~$C$ to get this strict inequality depends only on the function~$\mu_0$, by the intermediate of the functions $\delta (\e)$ and $\delta' (\e)$.
\hfill$\Box$


\section{The unique spreading speed case} \label{sec:unique}

We begin with some preliminary work that will be needed to estimate the spreading speeds. 
The proof of Theorem~\ref{thm-vitspread} is then separated into two parts: the first part (Section~\ref{vitspread-p1}) is devoted to the proof that $w^*\leq w_\infty$, while in the second part (Section~\ref{vitspread-p2}) we prove that $w_* \geq w_\infty$.

\subsection{Construction of the approximated eigenfunctions}\label{subsec:approxeigen}

For all $p\in\R$, we define
\begin{equation} \label{eq:defH} H(p):= \left\{ \begin{array}{ccl}
                                          j^{-1} (|p|) &\hbox{ if }& |p |\geq j(M),\\
M &\hbox{ if }& |p |< j(M).\\
\end{array}\right.\end{equation}
The fundamental property of this function is given by the following result. 

\begin{prop}\label{prop:exv} (Propositions $3.1$ and $3.2$ in \cite{HamelNadinRoques}) 
For all $p\in\R$, $H(p)$ is the unique real number such that there exists a continuous $1$-periodic viscosity solution~$v$ of 
\begin{equation}\label{HJ-eigen}
(v'(y)-p)^2+\mu_0 (y)= H(p) \hbox{ over } \R . 
\end{equation}
\end{prop}

Next, we will need, as a first step of our proof, the function $v$ given by Proposition \ref{prop:exv} to be piecewise $\mathcal{C}^2$. 
This is true under some non-degeneracy hypothesis on $\mu_0$.  We will check below in the second part of the proof of Theorem~\ref{thm-vitspread} that
it is always possible to assume that this hypothesis is satisfied by approximation. 

\begin{lem} \label{lem:regv}
Assume that $\mu_0 \in\mathcal{C}^2 (\R)$ and that 
\begin{equation} \label{hyp:nondegmu}
\hbox{if } \mu_0 (x_0)=\max_\R \mu_0, \hbox{ then } \mu_0''(x_0)<0. 
\end{equation}
Then for all $p\in\R$, equation (\ref{HJ-eigen}) admits a $1$-periodic solution $v_p\in W^{2,\infty}(\R)$ which is piecewise $\mathcal{C}^2 (\R)$. 
\end{lem}
{\bf Proof.}
The proof relies on the explicit formulation of $v_p$. Assume first that 
$p > j( M) = j\big(\|\mu_0 \|_\infty\big)$. Then it is easy to check (see \cite{HamelNadinRoques}) that 
\begin{equation} \label{explicitv} v_p(x):= p x -\int_0^x \sqrt{H(p)-\mu_0 (y)}dy\end{equation}
satisfies (\ref{HJ-eigen}). Then, the definition of $j$ implies that $v_p$ is $1$-periodic and, as $\mu_0\in \mathcal{C}^1 (\R)$ 
and $H(p) > \mu_0 (y)$ for all $y\in\R$, the function $v_p$ is $\mathcal{C}^2 (\R)$. 
The case $p < -j( M)$ is treated similarly. 

Next, if  $|p| \leq j( M)$, let $F$ define for all $Y\in [0,1]$ by:
$$F(Y):= p+\int_Y^1 \sqrt{M-\mu_0 (y)}dy -\int_0^Y \sqrt{M-\mu_0 (y)}dy.$$
Then $F$ is continuous and, as $|p| \leq j( M)$, 
$$F(0) =  p+\int_0^1 \sqrt{M-\mu_0 (y)}dy = p + j(M)\geq 0.$$
Similarly, $F(1) = p-j(M) \leq 0$. Thus, there exists $X\in [0,1]$ so that $F(X)=0$. We now define:
\begin{equation} v_p (x) = \left\{ 
\begin{array}{lcl} 
\ds px -\int_0^x \sqrt{M-\mu_0 (y)}dy &\hbox{ for all }& x\in [0,X], \\[0.3cm]
\ds px -\int_0^X \sqrt{M-\mu_0 (y)}dy +\int_X^x \sqrt{M-\mu_0 (y)}dy &\hbox{ for all }& x\in [X,1].
\end{array} \right. \end{equation}
From the definition of $X$, the function $v_p$ is $1$-periodic. It is continuous and derivable at any point $x\in [0,1)\backslash \{X\}$ with
$$ v_p'(x) = \left\{ \begin{array}{lcl} 
p-\sqrt{M-\mu_0 (x)} &\hbox{ for all }& x\in [0,X),\\
p+\sqrt{M-\mu_0 (x)} &\hbox{ for all }& x\in (X,1).\\
\end{array}\right. $$
Hence, it satisfies (\ref{HJ-eigen}) in the sense of viscosity solutions. 
Lastly, for all $x\in (0,X)$ so that $\mu_0 (x) \neq M$, one has
$$v_p''(x) = \frac{\mu_0'(x)}{2\sqrt{M -\mu_0 (x)}}.$$
If $\mu_0 (x_M)=M$, then (\ref{hyp:nondegmu}) implies that $\mu_0 (x) <M$ for all $x\neq x_M$ close to $x_M$ and a Taylor expansion gives
$$\lim_{x\to x_M, x\neq x_M}v_p''(x) = \sqrt{-\mu_0''(x_M)/2}.$$
Hence, $v_p''$ can be extended to a continuous function over $(0,X)$. Similarly, it can be extended over $(X,1)$. It follows that $v_p''$ is bounded over $[0,1]$
and that it is piecewise $\mathcal{C}^2 (\R)$. \hfill $\Box$

\bigskip

For any $p \in \R$, define the elliptic operator:
 $$L_p \varphi := \varphi'' - 2p\varphi ' + (p^2 + \mu_0 (\phi (x)))\varphi.  $$

\begin{lem}\label{lem-vepapprox}
For all $p\in\R$, let
\begin{equation}
 \varphi_p (x):= \exp \Big( \frac{v_p(\phi (x))}{\phi'(x)}\Big).
\end{equation}
If (\ref{hyp:nondegmu}) holds and $\mu_0 \in \C^2 (\R)$, then $\varphi_p$ is piecewise $\mathcal{C}^2 (\R)$ and one has 
\begin{equation}
 \frac{L_p \varphi_p(x) -H(p) \varphi_p(x)}{\varphi_p(x)} \to 0 \hbox{ as } x\to +\infty. 
\end{equation}
\end{lem}
{\bf Proof.} The function $\varphi_p$ is piecewise $\mathcal{C}^2 (\R)$ since $v_p$ is piecewise $\mathcal{C}^2 (\R)$. 
For all $x$ so that~$v_p$ is $\mathcal{C}^2$ in $x$, we can compute
\[
 \begin{array}{rcl}
  \ds\varphi_p'(x) & = & \lp\ds v'_p(\phi(x))-\frac{\phi''(x)}{(\phi'(x))^2}v_p(\phi(x))  \rp\varphi_p(x),\\
  \ds\varphi_p''(x)& = & \lp  \ds\phi'(x)v''_p(\phi(x))
                               - \ds\frac{\phi''(x)}{\phi'(x)}v'_p(\phi(x))  
                               + \ds\lp2\frac{(\phi''(x))^2}{(\phi'(x))^3}-\frac{\phi'''(x)}{(\phi'(x))^2}\rp v_p(\phi(x)) \rp \varphi_p(x)\\
                         &   & + \lp \ds v'_p(\phi(x))- \frac{\phi''(x)}{(\phi'(x))^2}v_p(\phi(x)) \rp^2 
                                \varphi_p(x).
 \end{array}
\]
This gives
\[
\begin{array}{rcl}
 \ds \frac{L_p\varphi_p(x)-H(p)\varphi_p(x)}{\varphi_p(x)} & = &
                                \ds\phi'(x)v''_p(\phi(x))
                               - \ds\frac{\phi''(x)}{\phi'(x)}v'_p(\phi(x)) 
                              + \ds\lp2\frac{(\phi''(x))^2}{(\phi'(x))^3}-\frac{\phi'''(x)}{(\phi'(x))^2}\rp v_p(\phi(x))\\[0.3cm]
                             && \ds-2\frac{\phi''(x)}{(\phi'(x))^2}v_p(\phi(x))v'_p(\phi(x))
                               + \lp\ds\frac{\phi''(x)}{(\phi'(x))^2}v_p(\phi(x))\rp^2+ v'_p(\phi(x))^2\\[0.3cm]
                            &&  \ds -2p \lp \ds  v'_p(\phi(x))-\frac{\phi''(x)}{(\phi'(x))^2}v_p(\phi(x))  \rp +p^2+\mu_0 (\phi(x))-H(p),\\[0.5cm]
                            &=&     \ds\phi'(x)v''_p(\phi(x))
                               - \ds\frac{\phi''(x)}{\phi'(x)}v'_p(\phi(x)) 
                              + \ds\lp2\frac{(\phi''(x))^2}{(\phi'(x))^3}-\frac{\phi'''(x)}{(\phi'(x))^2}\rp v_p(\phi(x))\\[0.3cm]
                             && \ds-2\frac{\phi''(x)}{(\phi'(x))^2}v_p(\phi(x))v'_p(\phi(x))
                               + \lp\ds\frac{\phi''(x)}{(\phi'(x))^2}v_p(\phi(x))\rp^2 \\[0.3cm]
		     && \ds+ 2p \frac{\phi''(x)}{(\phi'(x))^2}v_p(\phi(x)) .                
\end{array}
\]
As $v_p$ is periodic and $W^{2,\infty}$, $v''_p$ is bounded. It follows from \eqref{hyp:phi2} that 
\[
  \ds \frac{L_p\varphi_p(x)-H(p) \varphi_p(x)}{\varphi_p(x)} \to 0 \ \hbox{ as }  x\to+\infty.
\]
\hfill $\Box$

\begin{lem}\label{lem:cvv}
Define $\varphi_p$ as in Lemma \ref{lem-vepapprox}.
Then 
\begin{equation}\label{eq:cvv}
\frac{\ln \varphi_p (x)}{x} \to 0 \hbox{ as } x\to +\infty.
\end{equation}
\end{lem}
{\bf Proof.} 
One has 
\begin{equation}
\frac{\ln \varphi_p (x)}{x}=\frac{v_p (\phi (x))}{\phi'(x)x} \ \hbox{ for all }x\in\R.
\end{equation}
The function $x\mapsto v_p(\phi (x))$ is clearly bounded since $v_p$ is periodic. Hence, (\ref{hyp:phi2}) gives the conclusion. 
\hfill$\Box$


\subsection{Upper bound for the spreading speed}\label{vitspread-p1}

\noindent {\bf Proof of part 1 of Theorem \ref{thm-vitspread}.}  We first assume that $\mu_0 \in \C^2 (\R)$. Let us now show that $w^* \leq w_\infty$. Let $c>w_\infty$ and $c_1 \in (w_\infty,c)$. We know 
that there exists $p\geq j(M)>0$ such that 
$$w_\infty=\min_{p' \geq j(M) }H(p')/p'= H(p)/p.$$ 
Let $k\geq M$ so that $p=j(k)>0,$  and $\varphi_p$ defined as in Lemma~\ref{lem-vepapprox}.
We  know from Lemma~\ref{lem-vepapprox} that there exists $X>0$ such that:
\begin{equation}
|L_{p} \varphi_p(x) -k \varphi_p(x)| \leq (c_1-w_\infty) j(k)\varphi_p(x) \hbox{ for all } x>X.
\end{equation}
Let $\overline{u}$ be defined for all $(t,x) \in [0,+\infty) \times \R$ by:
$$\overline{u}(t,x):= \min \{1, \varphi_p (x)e^{-j(k) (x-h-c_1t)}\},$$
where $h\in\R$ is large enough so that $u_0 (x) \leq \overline{u}(0,x)$ for all $x\in \R $ (this is always possible since $u_0$ is compactly supported). 
Moreover, $\overline{u}(t,x)<1$ if and only if $\varphi_p (x)<e^{j(k) (x-h-c_1t)}$, which is equivalent to 
$ x-v_p (\phi (x))/(j(k)\phi '(x))> h+c_1t \geq h$. Lemma \ref{lem:cvv} yields that the left hand-side of this inequality goes to $+\infty$ as $x\to +\infty$. Hence, we
can always take $h$ large enough so that $\overline{u}(t,x) <1$ implies $x>X$.   
It follows that for all $(t,x) \in [0,+\infty) \times \R$ such that $\overline{u}(t,x)<1$, one has
$$ \begin{array}{rcl}
    \partial_t \overline{u} -\partial_{xx} \overline{u} -f(x,\overline{u}) &\geq & \partial_t \overline{u} -\partial_{xx} \overline{u} -\mu (x)\overline{u}\\[0.2cm]
&\geq & j(k) c_1 \overline{u} - L_{j(k)} \big(\varphi_p \big)(x) e^{-j(k) (x-h-c_1t)}\\[0.2cm]
&\geq & j(k) c_1 \overline{u} - k\overline{u} -(c_1-w_\infty) j(k)\overline{u}\\[0.2cm]
&\geq & \big( j(k) w_\infty - k\big)\overline{u} =0.
   \end{array}$$
It follows from the parabolic maximum principle that $\overline{u }(t,x) \geq u(t,x) $ for all $(t,x) \in [0,+\infty) \times \R$. 
Hence, for all given $x\in\R$, 
$$ u(t,x) \leq \varphi_p (x) e^{-j(k) (x-h-c_1t)} \ \hbox{ for all } t>0.$$
Let $\e>0$ so that $\e < j(k) ( c-c_1)/c$. Lemma~\ref{lem:cvv} yields that there exists $R>0$ so that for all $x>R, \ln \big(\varphi_p (x) \big) \leq \e x$. 
Let $T= R/c$ and take $t\geq T$ and $x\geq c t$. One has
$$\begin{array}{rcl} 
\ln \Big(\varphi_p (x) e^{-j(k) (x-h-c_1t)}\Big)&=& \ln \big(\varphi_p (x) \big) -j(k) (x-h-c_1t)\\[0.2cm]
&\leq & (\e-j(k)) x +j(k)(h+c_1 t)\\[0.2cm]
&\leq & \big(\e c+j(k)(c_1-c) \big) t +j(k)h\\[0.2cm]
   & \to & -\infty  \hbox{ as } t\to +\infty\\
  \end{array}$$
since $\e c<j(k)(c-c_1)$. 
Hence, 
$$\lim_{t\to +\infty} \max_{x\geq ct} \Big(\varphi_p (x) e^{-j(k) (x-h-c_1 t)}\Big)= 0 \hbox{ as } t\to +\infty,$$
which ends the proof in the case $\mu_0 \in \C^2 (\R)$.

Lastly, if $\mu_0\in\mathcal{C}^0(\R)$ is an arbitrary $1$-periodic function, then one easily concludes by smoothing $\mu_0$ from above. Indeed, one can find a sequence $(\mu_0^n )_n \in \C^2 (\R)^{\mathbb{N}}$ converging uniformly to $\mu_0$, and such that for all $n\in \mathbb{N}$ and $x \in \R$, $\mu_0 (x) \leq \mu_0^n (x)$.

It follows from the maximum principle that 
$$\lim_{t\to +\infty} \max_{x\geq ct}u(t,x)=0 \ \hbox{ for all } w > \min_{k\geq M} k/j^n(k),$$
where 
$ \ds j^n (k) = \int_0^1 \sqrt{ k-\mu_0^n (x)}dx \geq j (k) >0.$
Letting $n\to +\infty$, one gets
$$\lim_{t\to +\infty} \max_{x\geq ct}u(t,x)=0  \ \hbox{ for all } w > w_\infty,$$
which concludes the proof. \hfill $\Box$


\subsection{Lower bound on the spreading speed}\label{vitspread-p2}

\noindent {\bf Proof of part 2 of Theorem \ref{thm-vitspread}.} 
First, assume that $\mu_0\in\mathcal{C}^2(\R)$ satisfies (\ref{hyp:nondegmu}). Let $\varphi_p$ as in Lemma \ref{lem-vepapprox}. For all $\delta>0$, 
take $R$ large enough so that
$L_p \varphi_p \geq (H(p) -\delta) \varphi_p$ at any point of  $(R,+\infty)$ where $\varphi_p$ is piecewise $\mathcal{C}^2$. 
It is easy to derive from the proof of Lemma \ref{lem-vepapprox} that $\varphi_p'/\varphi_p$ is bounded and uniformly continuous. Take 
$C>0$ so that $|\varphi_p' (x)| \leq C \varphi_p (x)$ for all $x\in\R$. 

We need more regularity in order to apply the results of~\cite{BerestyckiNadin}. Consider a compactly supported nonnegative 
mollifier $\chi \in \mathcal{C}^\infty (\R)$ so that $\int_\R \chi = 1$ and define the convoled function
$\psi_p := \exp \big(\chi \star \ln\varphi_p\big) \in\mathcal{C}^2 (\R)$. 
One has 
$\psi_p'/\psi_p = \chi \star \Big(\varphi_p'/\varphi_p\Big). $
Hence, $|\psi_p' (x)| \leq C \psi_p (x)$ for all $x\in\R$ and, as $\varphi_p'/\varphi_p$ and $\big( \varphi_p'/\varphi_p \big)^2$ are uniformly continuous, 
up to some rescaling of $\chi$, we can assume that 
$$\Big\|\big|\chi \star \Big(\frac{\varphi_p'}{\varphi_p}\Big)\big|^2 -\big|\frac{\varphi_p'}{\varphi_p}\big|^2 \Big\|_\infty  \leq \delta, 
\ \Big\|\chi \star \Big( \big|\frac{\varphi_p'}{\varphi_p}\big|^2\Big)-\big|\frac{\varphi_p'}{\varphi_p}\big|^2\Big\|_\infty\leq \delta 
\hbox{ and } \|\chi\star \mu  -\mu \|_\infty \leq \delta.$$
We now compute
$$\frac{ \psi_p''}{\psi_p}=\Big|\frac{\psi_p'}{\psi_p}\Big|^2 +\chi \star \Big(\frac{\varphi_p''}{\varphi_p}-\Big|\frac{\varphi_p'}{\varphi_p}\Big|^2\Big) 
=\Big|\chi \star \Big(\frac{\varphi_p'}{\varphi_p}\Big)\Big|^2 -\chi \star \Big( \Big|\frac{\varphi_p'}{\varphi_p}\Big|^2\Big)+\chi \star \Big(\frac{\varphi_p''}{\varphi_p}\Big) 
\geq -2\delta +\chi \star \Big(\frac{\varphi_p''}{\varphi_p}\Big).$$
It follows that
$$ \begin{array}{rcl} \displaystyle \frac{L_p \psi_p}{\psi_p} = \frac{\psi_p''-2p \psi_p' + \mu (x)\psi_p}{\psi_p}&\geq& 
-2\delta+\chi \star \Big( \displaystyle\frac{\varphi_p''}{\varphi_p}-2p  \displaystyle\frac{\varphi_p'}{\varphi_p}\Big)+\mu (x)\\
&\geq & -2\delta+\chi \star \Big( H(p)-\mu-\delta \Big)+\mu (x)\\
&\geq & -3\delta + H(p)+\mu (x)-\chi \star \mu (x)\\
&\geq & -4\delta + H(p)\\
 \end{array}$$  
in $(R,+\infty)$.
On the other hand, Lemma~\ref{lem:cvv} yields $\psi_p \in\mathcal{A}_R$, where $\mathcal{A}_R$ is the set of admissible test-functions (in the sense of \cite{BerestyckiNadin}) over $(R,\infty)$:
\begin{equation} \label{defA} \begin{array}{rll}
 \A_R:= \big\{&\psi\in \mathcal{C}^0([R,\infty) )\cap\mathcal{C}^{2}((R,\infty) ),&\\ &\psi'/\psi \in L^\infty ((R,\infty)), \
 \psi>0 \hbox{ in } [R,\infty), \ \lim_{x\to +\infty} \frac{1}{x}\ln\psi (x) =0 & \big\}.\\ \end{array}
\end{equation}
Thus, one has $\underline{\lambda_1} (L_p ,(R,+\infty)) \geq H(p)-4\delta$, where the principal eigenvalue $\underline{\lambda_1}$ is defined by  
\begin{equation}\label{deflambda1'}
\underline{\lambda_1}(L_p,(R,\infty) ):=\sup\{\lambda\ |\
\exists\phi\in \A_R \hbox{ such that } L_p\phi\geq\lambda\phi \hbox{ in }(R,\infty) \},
\end{equation}
Hence, $\lim_{R\to +\infty} \underline{\lambda_1} (L_p ,(R,+\infty)) \geq H(p)$ for all $p>0$.

In order to use Theorem 2.1 of \cite{BerestyckiNadin}, we need the nonlinearity to have two steady states and to be positive between these two steady states. 
It is not the case here but we will bound $f$ from below by such a nonlinearity. 
As $\min_\R\mu_0>0$ and $f$ is of class $\mathcal{C}^1$ in the neighborhood of $u=0$, we know that there exists $\theta \in (0,1)$ so that 
$$f(x,u) >0 \hbox{ for all } x\in\R \ \hbox{ and } \ u\in (0,\theta).$$
Let $\zeta=\zeta (u)$ a smooth function so that 
$$0< \zeta (u)\leq 1 \hbox{ for all } u\in (0,\theta), \hspace{0.2cm} \zeta (u)=0 \hbox{ for all } u\geq \theta \ \hbox{ and } \ \zeta (u)=1\hbox{ for all }u\in (0,\frac{\theta}{2}).$$ 
Define $\underline{f} (x,u) := \zeta (u) f(x,u)$ for all $(x,u) \in \R\times [0,1]$. Then 
$$\underline{f} \leq f \hbox{ in } \R\times [0,1]\ \hbox{ and } \ \underline{f}_u' (x,0)= f_u'(x,0)=\mu_0 (\phi (x))\hbox{ for  all } x\in\R.$$ Let $\underline{u}$ the solution of (\ref{eqn:eqRD}) with nonlinearity $\underline{f}$ instead of $f$ and initial datum $u_0$. 
The parabolic maximum principle yields $u\geq \underline{u}$. 

Since the function $\underline{f}$ satisfies the hypotheses of Theorem 2.1 in \cite{BerestyckiNadin}, we conclude that
$$\lim_{t\to +\infty} \min_{x\in [0, wt]} \underline{u} (t,x) =1 \hbox{ for all } w\in \Big( 0, \min_{p>0} \displaystyle \frac{H(p)}{p}\Big).$$
It follows that 
$$ w_* \geq \min_{p>0} \displaystyle \frac{H(p)}{p} = \min_{k \geq M} \frac{k}{j(k)}.$$

Next, assume that $\mu_0\in\mathcal{C}^2(\R)$ does not satisfy (\ref{hyp:nondegmu}). 
Let $\overline{y}\in\R$ so that $\mu_0 (\overline{y}) = \ds\max_{y\in\R} \mu_0 (y)$. 
Take a $1$-periodic function $\chi \in\mathcal{C}^2 (\R)$ so that $\chi (0)=0$, $\chi (y) >0$ for all $y\neq 0$ and $\chi ''(0) >0$. Define for all $n\in\mathbb{N}$ and~$x\in\R$:
$$\mu_0^n (y) := \mu_0 (y) - \frac{1}{n} \chi (y-\overline{y}).$$
This $1$-periodic function satisfies (\ref{hyp:nondegmu}) for all $n$ and one has $0< \mu_0^n \leq \mu_0$ for~$n$ large enough. It follows from the maximum principle that 
$$\liminf_{t\to +\infty} \min_{0\leq x\leq wt}u(t,x) >0 \hbox{ for all } w\in \Big(0,\ds \min_{k\geq M}{ \frac{k}{j^n(k)}}\Big),$$
where 
$\ds j^n (k) = \int_0^1 \sqrt{ k-\mu_0^n (x)}dx \geq j (k) > 0$ for all $k \geq M$.
Letting $n\to +\infty$, one has $\mu_0^n (y) \to \mu_0 (y)$ uniformly in~$y\in\R$ and thus 
$$\liminf_{t\to +\infty} \min_{0\leq x\leq wt}u(t,x) >0 \hbox{ for all } w\in (0, w_\infty),$$
which concludes the proof in this case. 

\smallskip

Lastly, if $\mu_0\in\mathcal{C}^0(\R)$ is an arbitrary $1$-periodic function, then one easily concludes by smoothing $\mu_0$ as in the previous step. 
 \hfill $\Box$ 



\begin{thebibliography}{1}

\bibitem{AronsonWeinberger}
D.G. Aronson, and H.F. Weinberger.
\newblock Multidimensional nonlinear diffusions arising in population genetics.
\newblock {\em Adv. Math.}, 30:33--76, 1978.

\bibitem{BerestyckiHamel}
H.~Berestycki, and F.~Hamel.
\newblock Front propagation in periodic excitable media.
\newblock {\em Comm. Pure Appl. Math.}, 55:949--1032, 2002.

\bibitem{BHNa}
H.~Berestycki, F.~Hamel, and G.~Nadin.
\newblock Asymptotic spreading in heterogeneous diffusive excitable media.
\newblock {\em J. Func. Anal.}, 255(9):2146--2189, 2008.

\bibitem{BHR07}
H. ~Berestycki, F.~Hamel, and L.~Rossi.
\newblock Liouville type results for semilinear elliptic equations in unbounded domains.
\newblock {\em Annali Mat. Pura Appli.}, 186:469-507, 2007.

\bibitem{BerestyckiNadin}
H.~Berestycki, and G.~Nadin.
\newblock Spreading speeds for one-dimensional monostable reaction-diffusion equations.
\newblock {\em preprint}.

\bibitem{Freidlin}
M.~Freidlin.
\newblock On wave front propagation in periodic media.
\newblock {\em In: Stochastic analysis and applications, ed. M. Pinsky,
  Advances in Probability and related topics}, 7:147--166, 1984.

\bibitem{GartnerFreidlin}
M.~Freidlin, and J.~G\"artner.
\newblock On the propagation of concentration waves in periodic and random
  media.
\newblock {\em Sov. Math. Dokl.}, 20:1282--1286, 1979.

\bibitem{HamelNadin}
F. Hamel, and G. Nadin. 
\newblock Spreading properties and complex dynamics for monostable reaction-diffusion equations.
\newblock {\em preprint}.

\bibitem{HamelNadinRoques}
F.~Hamel, G.~Nadin and L. Roques.
\newblock A viscosity solution method for the spreading speed formula in slowly varying media.
\newblock {\em  Indiana Univ. Math. J}, to appear. 


\bibitem{MallordyRoquejoffre}
J.-F. Mallordy, and J.-M. Roquejoffre.
\newblock A parabolic equation of the KPP type in higher dimensions.
\newblock {\em SIAM J. Math. Anal.}, 26(1): 1--20, 1995. 

\bibitem{NolenRuddXin}
J.~Nolen, M.~Rudd, and J.~Xin.
\newblock Existence of KPP fronts in spatially-temporally periodic advection
  and variational principle for propagation speeds. 
\newblock {\em Dynamics of PDE}, 2(1):1--24, 2005.

\bibitem{NolenXinrandom}
J. Nolen, and J. Xin.
\newblock {Asymptotic Spreading of KPP Reactive Fronts in Incompressible Space-Time Random Flows}.
\newblock {\em Ann. de l'Inst. Henri Poincare -- Analyse Non Lineaire}, 26(3):815--839, 2008.

\bibitem{Weinberger}
H.~Weinberger.
\newblock On spreading speed and travelling waves for growth and migration
  models in a periodic habitat.
\newblock {\em J. Math. Biol.}, 45:511--548, 2002.

\end{thebibliography}
\end{document}